\documentclass[12pt]{article}%
\usepackage{amsfonts}
\usepackage{sw20bams}
\usepackage{amsmath}
\usepackage{amssymb}
\usepackage{graphicx}%
\providecommand{\U}[1]{\protect\rule{.1in}{.1in}}
\begin{document}

\title{Components and Cycles of Random Mappings}
\author{Steven Finch}
\date{May 11, 2022}
\maketitle

\begin{abstract}
Each connected component of a mapping $\{1,2,\ldots,n\}\rightarrow
\{1,2,\ldots,n\}$ contains a unique cycle. \ The largest such component can be
studied probabilistically via either a delay differential equation or an
inverse Laplace transform. \ The longest such cycle likewise admits two
approaches:\ we find an (apparently new) density formula for its length.
\ Implications of a constraint -- that exactly one component exists -- are
also examined. \ For instance, the mean length of the longest cycle is
$(0.7824...)\sqrt{n}$ in general, but for the special case, it is
$(0.7978...)\sqrt{n}$, a difference of less than$\ 2\%$.

\end{abstract}

\footnotetext{Copyright \copyright \ 2022 by Steven R. Finch. All rights
reserved.}Two delay differential equations (DDEs) shall be helpful:%
\[%
\begin{array}
[c]{ccc}%
x\,\rho^{\prime}(x)+\rho(x-1)=0\text{ for }x>1, &  & \rho(x)=1\text{ for
}0\leq x\leq1
\end{array}
\]
where $\rho$ is Dickman's function \cite{Dick-tcs0, More-tcs0, Lag-tcs0}, and%
\[%
\begin{array}
[c]{ccc}%
x\,\sigma^{\prime}(x)+\frac{1}{2}\sigma(x)+\frac{1}{2}\sigma(x-1)=0\text{ for
}x>1, &  & \sigma(x)=1/\sqrt{x}\text{ for }0<x\leq1
\end{array}
\]
where $\sigma$ could justifiably be called Watterson's function
\cite{Watt-tcs0, ABT1-tcs0, ABT2-tcs0}. \ It is understood that $\rho
(x)=0=\sigma(x)$ for $x<0$. \ One-sided Laplace transforms will also play a
role; for example, we have%
\[%
\begin{array}
[c]{ccccc}%
\mathcal{L}\left[  \rho(\xi)\right]  =\dfrac{\exp\left(  -E(\eta)\right)
}{\eta}, &  & \mathcal{L}\left[  \sigma(\xi)\right]  =\dfrac{\exp\left(
-\frac{1}{2}E(\eta)\right)  }{\sqrt{\eta/\pi}}, &  & \eta\in\mathbb{C}%
\smallsetminus(-\infty,0]
\end{array}
\]
or equivalently%
\[%
\begin{array}
[c]{ccc}%
\mathcal{L}^{-1}\left[  \dfrac{\exp\left(  -E(\eta)\right)  }{\eta}\right]
=\rho(\xi), &  & \mathcal{L}^{-1}\left[  \dfrac{\exp\left(  -\frac{1}{2}%
E(\eta)\right)  }{\sqrt{\eta/\pi}}\right]  =\sigma(\xi)
\end{array}
\]
where%
\[%
\begin{array}
[c]{ccc}%
E(x)=%
{\displaystyle\int\limits_{x}^{\infty}}
\dfrac{e^{-t}}{t}dt=-\operatorname{Ei}(-x), &  & x>0
\end{array}
\]
is the exponential integral \cite{Fi1-tcs0, Fi2-tcs0, Fi3-tcs0, Fi4-tcs0}.

\section{Stepanov}

For introductory purposes, let us examine solely one-to-one mappings
$\{1,2,\ldots,n\}\rightarrow\{1,2,\ldots,n\}$, i.e., permutations on $n$
symbols. \ Let $\Lambda$ denote the length of the longest cycle in an
$n$-permutation, chosen uniformly at random. \ We have%
\[%
\begin{array}
[c]{ccc}%
\lim\limits_{n\rightarrow\infty}\mathbb{P}\left\{  \Lambda\leq a\,n\right\}
=\rho\left(  \dfrac{1}{a}\right)  , &  & 0<a\leq1
\end{array}
\]
where, from the DDE,
\[
\rho(x)=\left\{
\begin{array}
[c]{lll}%
1-\ln(x) &  & \text{if }1<x\leq2,\\
1-\dfrac{\pi^{2}}{12}-\ln(x)+\dfrac{1}{2}\ln(x)^{2}+\operatorname*{Li}%
\nolimits_{2}\left(  \dfrac{1}{x}\right)  &  & \text{if }2<x\leq3.
\end{array}
\right.
\]
Thus, for example,
\begin{align*}
\lim\limits_{n\rightarrow\infty}\mathbb{P}\left\{  \frac{1}{3}<\frac{\Lambda
}{n}\leq\frac{1}{2}\right\}   &  =\rho(2)-\rho(3)\\
&  =\frac{\pi^{2}}{12}-\ln(2)+\ln(3)-\frac{1}{2}\ln(3)^{2}-\operatorname*{Li}%
\nolimits_{2}\left(  \frac{1}{3}\right) \\
&  =0.258244431148....
\end{align*}
Let us now explore a less-familiar approach \cite{Step-tcs0}. \ Define a
function $h(\xi)$ to be equal to $0$ for $\xi<1$ and $1/\xi$ for $\xi\geq1$.
\ Using the series expansion for $\exp\left(  -E(\eta)\right)  $ in terms of
$E(\eta)$, we have%
\[
\mathcal{L}^{-1}\left[  \dfrac{\exp\left(  -E(\eta)\right)  }{\eta}\right]  =%
{\displaystyle\sum\limits_{k=0}^{\infty}}
\frac{(-1)^{k}}{k!}\mathcal{L}^{-1}\left[  \dfrac{E(\eta)^{k}}{\eta}\right]
\]
where the power $E(\eta)^{k}$ is the Laplace transform of the $k$-fold
self-convolution $h_{k}(\xi)$ of $h(\xi)$. \ On the one hand, formulas%
\[%
\begin{array}
[c]{ccc}%
\mathcal{L}^{-1}\left[  \dfrac{1}{\eta}\right]  =1, &  & \mathcal{L}%
^{-1}\left[  \dfrac{E(\eta)}{\eta}\right]  =\left\{
\begin{array}
[c]{lll}%
\ln(\xi) &  & \text{if }\xi\geq1,\\
0 &  & \text{if }0\leq\xi<1
\end{array}
\right.
\end{array}
\]
are well-known. \ On the other hand,%
\[
\mathcal{L}^{-1}\left[  \dfrac{E(\eta)^{2}}{\eta}\right]  =\left\{
\begin{array}
[c]{lll}%
-\dfrac{\pi^{2}}{6}+\ln(\xi)^{2}+2\operatorname*{Li}\nolimits_{2}\left(
\dfrac{1}{\xi}\right)  &  & \text{if }\xi\geq2,\\
0 &  & \text{if }0\leq\xi<2.
\end{array}
\right.
\]
does not appear in \cite{OB-tcs0}. \ Since $h_{k}(\xi)=0$ for $0\leq\xi<k$,
only terms $k=0,1,\ldots,\left\lfloor \xi\right\rfloor $ of the series need to
be summed. \ Consequently, for $1/3<a\leq1/2$,
\[
\mathcal{L}^{-1}\left[  \dfrac{\exp\left(  -E(\eta)\right)  }{\eta}\right]
_{\xi\rightarrow1/a}=1-\ln\left(  \frac{1}{a}\right)  +\frac{1}{2}\left[
-\dfrac{\pi^{2}}{6}+\ln(a)^{2}+2\operatorname*{Li}\nolimits_{2}\left(
a\right)  \right]
\]
which is consistent with before.

\section{Mutafchiev}

Let us now remove the restriction that mappings be one-to-one. \ Let $\Lambda$
denote the length of the largest component in an $n$-mapping, chosen uniformly
at random. \ We have
\[%
\begin{array}
[c]{ccc}%
\lim\limits_{n\rightarrow\infty}\mathbb{P}\left\{  \Lambda\leq a\,n\right\}
=\dfrac{1}{\sqrt{a}}\;\sigma\left(  \dfrac{1}{a}\right)  , &  & 0<a\leq1
\end{array}
\]
where, from the DDE,%
\[%
\begin{array}
[c]{ccc}%
\sqrt{x}\,\sigma(x)=1-\dfrac{1}{2}\ln\left(  \dfrac{1+\sqrt{1-\frac{1}{x}}%
}{1-\sqrt{1-\frac{1}{x}}}\right)  &  & \text{if }1<x\leq2.
\end{array}
\]
For reasons of simplicity, change our example domain from $[1/3,1/2]$ to
$[1/2,2/3]$. \ Hence
\begin{align*}
\lim\limits_{n\rightarrow\infty}\mathbb{P}\left\{  \frac{1}{2}<\frac{\Lambda
}{n}\leq\frac{2}{3}\right\}   &  =\sqrt{\frac{3}{2}}\;\sigma\left(  \frac
{3}{2}\right)  -\sqrt{2}\;\sigma(2)\\
&  =-\frac{1}{2}\ln\left(  \dfrac{1+\sqrt{1-\frac{2}{3}}}{1-\sqrt{1-\frac
{2}{3}}}\right)  +\frac{1}{2}\ln\left(  \dfrac{1+\sqrt{1-\frac{1}{2}}}%
{1-\sqrt{1-\frac{1}{2}}}\right) \\
&  =0.222894638557....
\end{align*}
Let us again explore the less-familiar approach \cite{Muta-tcs0}. \ Define
$h(\xi)$ as previously. \ From%
\[
\sqrt{\pi\,\xi}\;\mathcal{L}^{-1}\left[  \dfrac{\exp\left(  -\frac{1}{2}%
E(\eta)\right)  }{\sqrt{\eta}}\right]  =\sqrt{\pi\,\xi}\;%
{\displaystyle\sum\limits_{k=0}^{\infty}}
\frac{(-1)^{k}}{2^{k}k!}\mathcal{L}^{-1}\left[  \dfrac{E(\eta)^{k}}{\sqrt
{\eta}}\right]
\]
we recognize two well-known formulas:
\[%
\begin{array}
[c]{ccc}%
\mathcal{L}^{-1}\left[  \dfrac{1}{\sqrt{\eta}}\right]  =\dfrac{1}{\sqrt
{\pi\,\xi}}, &  & \mathcal{L}^{-1}\left[  \dfrac{E(\eta)}{\sqrt{\eta}}\right]
=\left\{
\begin{array}
[c]{lll}%
\dfrac{2}{\sqrt{\pi\,\xi}}\operatorname{arctanh}\left(  \sqrt{1-\dfrac{1}{\xi
}}\right)  &  & \text{if }\xi\geq1,\\
0 &  & \text{if }0\leq\xi<1.
\end{array}
\right.
\end{array}
\]
Since $h_{k}(\xi)=0$ for $0\leq\xi<k$, only terms $k=0,1,\ldots,\left\lfloor
\xi\right\rfloor $ of the series need to be summed. \ Consequently, for
$1/2<a\leq1$,%
\[
\sqrt{\frac{\pi}{a}}\;\mathcal{L}^{-1}\left[  \dfrac{\exp\left(  -\frac{1}%
{2}E(\eta)\right)  }{\sqrt{\eta}}\right]  _{\xi\rightarrow1/a}%
=1-\operatorname{arctanh}\left(  \sqrt{1-a}\right)  =1-\frac{1}{2}\ln\left(
\dfrac{1+\sqrt{1-a}}{1-\sqrt{1-a}}\right)
\]
which again is consistent with before.

As an aside, had we kept the domain as $[1/3,1/2]$, then numerics are
possible:
\begin{align*}
\lim\limits_{n\rightarrow\infty}\mathbb{P}\left\{  \frac{1}{3}<\frac{\Lambda
}{n}\leq\frac{1}{2}\right\}   &  =%
{\displaystyle\int\limits_{1/3}^{1/2}}
\dfrac{1}{2\,t^{3/2}}\;\sigma\left(  \dfrac{1-t}{t}\right)  dt\\
&  =\frac{1}{2}%
{\displaystyle\int\limits_{1/3}^{1/2}}
\dfrac{1}{t\sqrt{1-t}}\,\left[  1-\dfrac{1}{2}\ln\left(  \dfrac{1+\sqrt
{1-\frac{t}{t-1}}}{1-\sqrt{1-\frac{t}{t-1}}}\right)  \right]  dt\\
&  =0.110414874191....
\end{align*}
but symbolics seem unlikely. \ A closed-form expression for the inverse
Laplace transform of $E(\eta)^{2}/\sqrt{\eta}$ also remains open.

\section{Purdom \& Williams}

We change the subject to cycles, i.e., loops in the functional graph. Let
$\Lambda$ denote the length of the longest cycle in an $n$-mapping, chosen
uniformly at random. \ Our goal is to find $\lim\nolimits_{n\rightarrow\infty
}\mathbb{P}\left\{  \Lambda\leq\,b\sqrt{n}\right\}  $. \ While no relevant DDE
is yet known, there is an associated inverse Laplace transform%
\[
\sqrt{\frac{\pi}{b}}\;\mathcal{L}^{-1}\left[  \dfrac{\exp\left(  -E\left(
\sqrt{2\,b\,\eta}\right)  \right)  }{\sqrt{\eta}}\right]  _{\xi\rightarrow
1/b}=\sqrt{\frac{\pi}{b}}\;%
{\displaystyle\sum\limits_{k=0}^{\infty}}
\frac{(-1)^{k}}{k!}\mathcal{L}^{-1}\left[  \dfrac{E\left(  \sqrt{2\,b\,\eta
}\right)  ^{k}}{\sqrt{\eta}}\right]  _{\xi\rightarrow1/b}%
\]
due to Mutafchiev \cite{Muta-tcs0}. \ Unlike previously, $0<b<\infty$
holds\ (rather than $0<a<1$) and $E\left(  \sqrt{2\,b\,\eta}\right)  $ is the
Laplace transform of $\frac{1}{2\,\xi}\operatorname*{erfc}\left(  \sqrt
{\frac{b}{2\,\xi}}\right)  $. \ Self-convolutions of this function do not
enjoy the same vanishing properties as those for $h(\xi)$. \ Truncating the
infinite series, although it is convergent, will unfortunately lead to
non-zero error. \ An alternative expression
\[
\sqrt{2\pi}\,\frac{1}{2\pi i}%
{\displaystyle\int\limits_{1-i\infty}^{1+i\infty}}
\exp\left(  -E(b\,\zeta)+\frac{\zeta^{2}}{2}\right)  d\zeta
\]
is due to Stepanov \cite{Step-tcs0}. \ Letting $\eta=\frac{b\,\zeta^{2}}{2}$,
i.e., $\zeta=\sqrt{\frac{2\,\eta}{b}}$ and $d\zeta=\frac{d\eta}{\sqrt
{2\,b\,\eta}}$, demonstrates that this complex contour integral is equal to
the preceding inverse Laplace transform.

There is, thankfully, a different approach available. \ If the number $N$\ of
cyclic points in a random mapping is fixed, then as Kolchin \cite{Kolc-tcs0}
wrote, \textquotedblleft... these $N$ cyclic points... form a random
permutation\textquotedblright. \ This suggests multiplying \cite{Pap-tcs0} the
conditional probability of $\Lambda$ given $N$:%
\[
\mathbb{P}\left\{  \Lambda\leq\lambda\sqrt{n}\;|\;N=\nu\sqrt{n}\right\}
=\rho\left(  \frac{\nu}{\lambda}\right)
\]
by \cite{RS-tcs0, Har-tcs0} the limiting density (as $n\rightarrow\infty$) of
$N$:%
\[%
\begin{array}
[c]{ccc}%
\nu\exp\left(  \dfrac{-\nu^{2}}{2}\right)  & \;\;\;\;\;\;\;\;\;\; &
\text{(Rayleigh)}%
\end{array}
\]
and then differentiating (with respect to $\lambda$) to obtain the joint
density of $(\Lambda,N)$:%
\begin{align*}
f(\lambda,\nu)  &  =\nu\exp\left(  \frac{-\nu^{2}}{2}\right)  \rho^{\prime
}\left(  \frac{\nu}{\lambda}\right)  \left(  \frac{-\nu}{\lambda^{2}}\right)
\\
&  =\nu\exp\left(  \frac{-\nu^{2}}{2}\right)  \rho\left(  \frac{\nu}{\lambda
}-1\right)  \left(  \frac{-\lambda}{\nu}\right)  \left(  \frac{-\nu}%
{\lambda^{2}}\right) \\
&  =\frac{\nu}{\lambda}\exp\left(  \frac{-\nu^{2}}{2}\right)  \rho\left(
\frac{\nu-\lambda}{\lambda}\right)
\end{align*}
where $0<\lambda<\nu<\infty$. \ We have not seen this formula in the
literature:\ it is apparently new. \ For $r\geq2$, let $\Lambda_{r}$ denote
the length of the $r^{\text{th}}$ longest cycle in an $n$-mapping, chosen
uniformly at random. \ If the permutation has no $r^{\text{th}}$ cycle, then
its $r^{\text{th}}$ longest cycle is defined to have length $0$. \ Define the
$r^{\text{th}}$ generalized Dickman function $\rho_{r}(\xi)$ to satisfy
\[%
\begin{array}
[c]{ccc}%
\xi\,\rho_{r}^{\prime}(\xi)+\rho_{r}(\xi-1)=\rho_{r-1}(\xi-1)\text{ for }%
\xi>1, &  & \rho_{r}(\xi)=1\text{ for }0\leq\xi\leq1
\end{array}
\]
where $\rho_{1}=\rho$. \ By the same argument, for $r\geq2$,%
\[
f_{r}(\lambda,\nu)=\frac{\nu}{\lambda}\exp\left(  \frac{-\nu^{2}}{2}\right)
\left[  \rho_{r}\left(  \frac{\nu-\lambda}{\lambda}\right)  -\rho_{r-1}\left(
\frac{\nu-\lambda}{\lambda}\right)  \right]  .
\]
Define%
\[
G_{r,h}=\frac{1}{h!(r-1)!}%
{\displaystyle\int\limits_{0}^{\infty}}
x^{h-1}E(x)^{r-1}\exp\left[  -E(x)-x\right]  dx
\]
(in this paper, rank $r=1,2,3$ or $4$; height $h=1$ or $2$). \ Purdom
\&\ Williams \cite{PW-tcs0}, building on Shepp \&\ Lloyd \cite{SL-tcs0},
discovered asymptotic formulas for moments $\mathbb{E}\left(  \Lambda_{r}%
^{h}\right)  $. \ We can easily verify their findings:%
\[
\lim_{n\rightarrow\infty}\frac{\mathbb{E}\left(  \Lambda_{r}\right)  }%
{\sqrt{n}}=\sqrt{\frac{\pi}{2}}\;G_{r,1}=\left\{
\begin{array}
[c]{ccc}%
0.78248160099165661501... &  & \text{if }r=1,\\
0.26267067265131265469... &  & \text{if }r=2,\\
0.11068781528281010827... &  & \text{if }r=3,\\
0.05056118481134243184... &  & \text{if }r=4;
\end{array}
\right.
\]%
\[
\lim_{n\rightarrow\infty}\frac{\mathbb{V}\left(  \Lambda_{r}\right)  }%
{n}=2G_{r,2}-\frac{\pi}{2}G_{r,1}^{2}=\left\{
\begin{array}
[c]{ccc}%
0.24111407342881901748... &  & \text{if }r=1,\\
0.04395998473216610374... &  & \text{if }r=2,\\
0.01233552055537805858... &  & \text{if }r=3,\\
0.00386619224804518754... &  & \text{if }r=4.
\end{array}
\right.
\]
The mode of $\Lambda_{1}$ occurs when%
\begin{align*}
0  &  =\frac{d}{d\lambda}%
{\displaystyle\int\limits_{\lambda}^{\infty}}
f(\lambda,\nu)d\nu=-f(\lambda,\lambda)+%
{\displaystyle\int\limits_{\lambda}^{\infty}}
\frac{\partial f}{\partial\lambda}(\lambda,\nu)d\nu\\
&  =-\exp\left(  \frac{-\lambda^{2}}{2}\right)  +%
{\displaystyle\int\limits_{\lambda}^{\infty}}
\nu\exp\left(  \frac{-\nu^{2}}{2}\right)  \frac{\partial}{\partial\lambda
}\left[  \frac{1}{\lambda}\rho\left(  \frac{\nu}{\lambda}-1\right)  \right]
d\nu;
\end{align*}
the inner derivative becomes%
\begin{align*}
\frac{-1}{\lambda^{2}}\rho\left(  \frac{\nu}{\lambda}-1\right)  +\frac
{1}{\lambda}\rho^{\prime}\left(  \frac{\nu}{\lambda}-1\right)  \left(
\frac{-\nu}{\lambda^{2}}\right)   &  =\frac{-1}{\lambda^{2}}\rho\left(
\frac{\nu}{\lambda}-1\right)  -\frac{1}{\lambda}\frac{\rho\left(  \frac{\nu
}{\lambda}-2\right)  }{\frac{\nu}{\lambda}-1}\left(  \frac{-\nu}{\lambda^{2}%
}\right) \\
&  =\frac{-1}{\lambda^{2}}\rho\left(  \frac{\nu-\lambda}{\lambda}\right)
+\frac{\nu}{\lambda^{2}(\nu-\lambda)}\rho\left(  \frac{\nu-2\lambda}{\lambda
}\right)  ;
\end{align*}
solving the equation%
\[
\exp\left(  \frac{-\lambda^{2}}{2}\right)  =\frac{1}{\lambda^{2}}%
{\displaystyle\int\limits_{\lambda}^{\infty}}
\left[  -\rho\left(  \frac{\nu-\lambda}{\lambda}\right)  +\frac{\nu}%
{\nu-\lambda}\rho\left(  \frac{\nu-2\lambda}{\lambda}\right)  \right]  \nu
\exp\left(  \frac{-\nu^{2}}{2}\right)  d\nu
\]
yields $0.4809...$ as the mode. \ The median $0.6842...$ of $\Lambda_{1}$
arises simply from%
\[
\frac{1}{2}=%
{\displaystyle\int\limits_{\lambda}^{\infty}}
{\displaystyle\int\limits_{\mu}^{\infty}}
f(\mu,\nu)d\nu\,d\mu.
\]
For $\Lambda_{2}$, the mode is $0$; we did not pursue the median. Another new
asymptotic result is the cross-correlation between $\Lambda_{r}$ and $N$:%
\[
\lim_{n\rightarrow\infty}\frac{\mathbb{E}\left(  \Lambda_{r}N\right)
-\mathbb{E}\left(  \Lambda_{r}\right)  \mathbb{E}\left(  N\right)  }%
{\sqrt{\mathbb{V}\left(  \Lambda_{r}\right)  }\sqrt{\mathbb{V}\left(
N\right)  }}=\frac{\sqrt{2-\frac{\pi}{2}}\;G_{r,1}}{\sqrt{2G_{r,2}-\frac{\pi
}{2}G_{r,1}^{2}}}=\left\{
\begin{array}
[c]{ccc}%
0.83298010... &  & \text{if }r=1,\\
0.65486924... &  & \text{if }r=2,\\
0.52094617... &  & \text{if }r=3,\\
0.42505712... &  & \text{if }r=4.
\end{array}
\right.
\]
Using formulas in \cite{Fi4-tcs0, Grf-tcs0, Shi-tcs0}, it is possible to
similarly compute the cross-correlation between $\Lambda_{r}$ and $\Lambda
_{s}$ where $r<s$.

\section{R\'{e}nyi}

A mapping is said to be\textbf{ connected} (or \textbf{indecomposable}) if it
possesses exactly one component. \ This is a rare event, in the sense that%
\[%
\begin{array}
[c]{ccc}%
\mathbb{P}\left\{  M=1\right\}  \sim\sqrt{\dfrac{\pi}{2n}} &  & \text{as
}n\rightarrow\infty
\end{array}
\]
where $M$ counts the components. \ Let $\Lambda$ denote the length of the
unique cycle in a connected $n$-mapping, chosen uniformly at random. \ Our
goal is to find $\lim\nolimits_{n\rightarrow\infty}\mathbb{P}\left\{
\Lambda\leq\,b\sqrt{n}\right\}  $ as before, but the circumstances are vastly
simpler. \ R\'{e}nyi \cite{Ren-tcs0} proved that the limiting density (as
$n\rightarrow\infty$) of $\Lambda$\ is%
\[%
\begin{array}
[c]{ccc}%
\sqrt{\dfrac{2}{\pi}}\exp\left(  \dfrac{-\lambda^{2}}{2}\right)  &
\;\;\;\;\;\;\;\;\;\; & \text{(half-normal)}%
\end{array}
\]
for $0<\lambda<\infty$, which implies immediately that%
\[
\lim_{n\rightarrow\infty}\frac{\mathbb{E}\left(  \Lambda\right)  }{\sqrt{n}%
}=\sqrt{\frac{2}{\pi}}=0.79788456080286535587...,
\]%
\[
\lim_{n\rightarrow\infty}\frac{\mathbb{V}\left(  \Lambda\right)  }{n}%
=1-\frac{2}{\pi}=0.36338022763241865692....
\]
It is surprising that arbitrary mappings and connected mappings differ so
little here ($0.7824...$ versus $0.7978...$). \ We might have expected that
uniqueness would carry more influence. \ 

Of course, $N=\Lambda$ when there is just one component. \ Allowing instead
$m$ components, where $m\geq2$ is a fixed integer, rarity persists
\cite{O1-tcs0, Pav1-tcs0}%
\[%
\begin{array}
[c]{ccc}%
\mathbb{P}\left\{  M=m\right\}  \sim\dfrac{1}{2^{m-1}(m-1)!}\sqrt{\dfrac{\pi
}{2n}}\ln(n)^{m-1} &  & \text{as }n\rightarrow\infty
\end{array}
\]
but the asymptotic values $\sqrt{2/\pi}$ and $1-2/\pi$ for $\mathbb{E}\left(
N\right)  $ and $\mathbb{V}\left(  N\right)  $ evidently do not change.
\ Section 6 contains more on this issue.

\section{Pavlov}

For arbitrary mappings, the expected number of components \cite{Kru-tcs0,
Rio-tcs0} is $\sim\frac{1}{2}\ln(n)$. \ If our constraint from Section 4
loosens so that $m\rightarrow\infty$ but so slowly that $m/\ln(n)\rightarrow
0$, then R\'{e}nyi's formula still applies, as proved by Pavlov
\cite{Pav1-tcs0}. \ This leads us to a set of conjectural results comparable
to those in Section 3.

Let $\Lambda$ be the longest cycle length and $N$ be the cyclic points total.
\ As before, a conditional probability coupled with the limiting density (as
$n\rightarrow\infty$) of $N$:
\[
\sqrt{\frac{2}{\pi}}\exp\left(  \frac{-\nu^{2}}{2}\right)
\]
suffice to give the joint density of $(\Lambda,N)$:%
\[
f(\lambda,\nu)=\sqrt{\frac{2}{\pi}}\;\frac{1}{\lambda}\exp\left(  \frac
{-\nu^{2}}{2}\right)  \rho\left(  \frac{\nu-\lambda}{\lambda}\right)
\]
where $0<\lambda<\nu<\infty$. \ Further,%
\[
f_{r}(\lambda,\nu)=\sqrt{\frac{2}{\pi}}\;\frac{1}{\lambda}\exp\left(
\frac{-\nu^{2}}{2}\right)  \left[  \rho_{r}\left(  \frac{\nu-\lambda}{\lambda
}\right)  -\rho_{r-1}\left(  \frac{\nu-\lambda}{\lambda}\right)  \right]
\]
for $r\geq2$. \ Moments are%
\[
\lim_{n\rightarrow\infty}\frac{\mathbb{E}\left(  \Lambda_{r}\right)  }%
{\sqrt{n}}=\sqrt{\frac{2}{\pi}}\;G_{r,1}=\left\{
\begin{array}
[c]{ccc}%
0.49814325870512904597... &  & \text{if }r=1,\\
0.16722134383091813637... &  & \text{if }r=2,\\
0.07046605176920746245... &  & \text{if }r=3,\\
0.03218824996523203019... &  & \text{if }r=4;
\end{array}
\right.
\]%
\[
\lim_{n\rightarrow\infty}\frac{\mathbb{V}\left(  \Lambda_{r}\right)  }%
{n}=G_{r,2}-\frac{2}{\pi}G_{r,1}^{2}=\left\{
\begin{array}
[c]{ccc}%
0.17854905846627743895... &  & \text{if }r=1,\\
0.02851495566901143371... &  & \text{if }r=2,\\
0.00732819205178914862... &  & \text{if }r=3,\\
0.00217522939296169629... &  & \text{if }r=4.
\end{array}
\right.
\]
The mode of $\Lambda_{1}$ occurs at $0$; the median at $0.3903...$. \ For
$\Lambda_{2}$, we did not pursue the median. \ The cross-correlation between
$\Lambda_{r}$ and $N$ is%
\[
\lim_{n\rightarrow\infty}\frac{\mathbb{E}\left(  \Lambda_{r}N\right)
-\mathbb{E}\left(  \Lambda_{r}\right)  \mathbb{E}\left(  N\right)  }%
{\sqrt{\mathbb{V}\left(  \Lambda_{r}\right)  }\sqrt{\mathbb{V}\left(
N\right)  }}=\frac{\sqrt{1-\frac{2}{\pi}}\;G_{r,1}}{\sqrt{G_{r,2}-\frac{2}%
{\pi}G_{r,1}^{2}}}=\left\{
\begin{array}
[c]{ccc}%
0.89066843... &  & \text{if }r=1,\\
0.74816251... &  & \text{if }r=2,\\
0.62190221... &  & \text{if }r=3,\\
0.52141727... &  & \text{if }r=4
\end{array}
\right.
\]
and, again, it is possible to compute the cross-correlation between
$\Lambda_{r}$ and $\Lambda_{s}$ where $r<s$.

The mean $0.4981...$ is sharply less than the other means $0.7824...$ and
$0.7978...$ we have exhibited. \ Why should this counterintuitive fact be
true? \ The scenario $m/\ln(n)\rightarrow0$ is intermediate to the others.
\ This is why we describe our work here as conjectural.

If instead $m/\ln(n)\rightarrow c$ for some constant $0<c<\infty$, then
Pavlov's \cite{Pav2-tcs0} density formula is%
\[
\frac{2^{c}\,\Gamma(c)}{\sqrt{2\pi}\,\Gamma(2c)}\nu^{2c}\exp\left(  \frac
{-\nu^{2}}{2}\right)  .
\]
This reduces to the density found in \cite{RS-tcs0, Har-tcs0} when $c=1/2$.

Given an arbitrary mapping, the \textbf{deepest} cycle is contained within the
largest component, whereas the \textbf{richest} component contains the longest
cycle. \ The deepest cycle need not be longest; the richest component need not
be largest. \ What can be said about the probability of either event, or the
average size of either structure? \ Questions about interplay at this level
appear to be difficult to answer.

For completeness, we mention \cite{De1-tcs0, De2-tcs0, RZ1-tcs0, RZ2-tcs0,
MPQS-tcs0}, which may offer additional insights and paths forward.

\section{Flajolet \& Odlyzko}

Let $a_{nm\ell}$ denote the number of $n$-mappings possessing exactly $m$
components and exactly $\ell$ cyclic points, where $n\geq\ell\geq m\geq2$.
\ We have \cite{Kup-tcs0, Eng-tcs0, O2-tcs0}%
\[%
{\displaystyle\sum\limits_{n=\ell}^{\infty}}
\,%
{\displaystyle\sum\limits_{\ell=m}^{n}}
\frac{a_{nm\ell}}{n!}x^{n}y^{\ell}=\frac{1}{m!}\ln\left(  \frac{1}%
{1-y\,\tau(x)}\right)  ^{m}%
\]
where $\tau(x)=x\exp(\tau(x))$ is Cayley's tree function. \ The dominant
singularity of $\tau(x)$ is at $x=e^{-1}$ and%
\[%
\begin{array}
[c]{ccc}%
\tau(x)\sim1-2^{1/2}\sqrt{1-e\,x} &  & \text{as }x\rightarrow e^{-1}.
\end{array}
\]
Differentiating with respect to $y$:%
\[%
{\displaystyle\sum\limits_{n=\ell}^{\infty}}
\,%
{\displaystyle\sum\limits_{\ell=m}^{n}}
\frac{\ell\,a_{nm\ell}}{n!}x^{n}y^{\ell-1}=\frac{1}{(m-1)!}\frac{\tau
(x)}{1-y\,\tau(x)}\ln\left(  \frac{1}{1-y\,\tau(x)}\right)  ^{m-1}%
\]
and setting $y=1$:%
\begin{align*}
\frac{1}{(m-1)!}\frac{\tau(x)}{1-\tau(x)}\ln\left(  \frac{1}{1-\tau
(x)}\right)  ^{m-1}  &  \sim\frac{1}{(m-1)!}\frac{1}{2^{1/2}\sqrt{1-e\,x}}%
\ln\left(  \frac{1}{2^{1/2}\sqrt{1-e\,x}}\right)  ^{m-1}\\
&  \sim\frac{1}{2^{m-1/2}(m-1)!}\sqrt{\frac{1}{1-e\,x}}\ln\left(  \frac
{1}{1-e\,x}\right)  ^{m-1}%
\end{align*}
we deduce%
\begin{align*}%
{\displaystyle\sum\limits_{\ell=m}^{n}}
\frac{\ell\,a_{nm\ell}}{n!}  &  \sim\frac{1}{2^{m-1/2}(m-1)!}\left(  \frac
{1}{e}\right)  ^{-n}\frac{\sqrt{n}\ln(n)^{m-1}}{n\,\Gamma(1/2)}\\
&  \sim\frac{\ln(n)^{m-1}}{2^{m-1}(m-1)!}\frac{e^{n}}{\sqrt{2\,\pi\,n}}%
\end{align*}
by the singularity analysis theorem of Flajolet \& Odlyzko \cite{FO-tcs0}.
\ Multiplying both sides by $n!/n^{n}$ and using Stirling's approximation, we
obtain%
\[%
{\displaystyle\sum\limits_{\ell=m}^{n}}
\frac{\ell\,a_{nm\ell}}{n^{n}}\sim\frac{\ln(n)^{m-1}}{2^{m-1}(m-1)!}.
\]
From Section 4,%
\[%
{\displaystyle\sum\limits_{\ell=m}^{n}}
\frac{a_{nm\ell}}{n^{n}}\sim\dfrac{\ln(n)^{m-1}}{2^{m-1}(m-1)!}\sqrt
{\dfrac{\pi}{2n}}%
\]
and hence, forming the ratio, $\mathbb{E}\left(  N\right)  \sim\sqrt{2n/\pi}$
as $n\rightarrow\infty$.

Differentiating again and setting $y=1$:%
\begin{align*}%
{\displaystyle\sum\limits_{n=\ell}^{\infty}}
\,%
{\displaystyle\sum\limits_{\ell=m}^{n}}
\frac{\ell(\ell-1)a_{nm\ell}}{n!}x^{n}  &  \sim\frac{1}{(m-1)!}\left(
\frac{\tau(x)}{1-\tau(x)}\right)  ^{2}\ln\left(  \frac{1}{1-\tau(x)}\right)
^{m-1}\\
&  \sim\frac{1}{(m-1)!}\frac{1}{2(1-e\,x)}\ln\left(  \frac{1}{2^{1/2}%
\sqrt{1-e\,x}}\right)  ^{m-1}\\
&  \sim\frac{1}{2^{m}(m-1)!}\frac{1}{1-e\,x}\ln\left(  \frac{1}{1-e\,x}%
\right)  ^{m-1}%
\end{align*}
we deduce%
\begin{align*}%
{\displaystyle\sum\limits_{\ell=m}^{n}}
\frac{\ell(\ell-1)a_{nm\ell}}{n!}  &  \sim\frac{1}{2^{m}(m-1)!}\left(
\frac{1}{e}\right)  ^{-n}\frac{n\ln(n)^{m-1}}{n\,\Gamma(1)}\\
&  \sim\frac{\ln(n)^{m-1}}{2^{m}(m-1)!}e^{n}.
\end{align*}
Multiplying by $n!/n^{n}$ and via the preceding, we obtain%
\[%
{\displaystyle\sum\limits_{\ell=m}^{n}}
\frac{\ell^{2}a_{nm\ell}}{n^{n}}\sim\frac{\ln(n)^{m-1}}{2^{m}(m-1)!}%
\sqrt{2\,\pi\,n}\sim\frac{\ln(n)^{m-1}}{2^{m-1}(m-1)!}\sqrt{\frac{\pi\,n}{2}%
}.
\]
Forming the ratio, $\mathbb{E}\left(  N^{2}\right)  \sim n$ as $n\rightarrow
\infty$ and thus $\mathbb{V}\left(  N\right)  \sim(1-2/\pi)n$.

\section{Addendum:\ Divisibility}

Let $m$ be a positive integer. \ A random variable $X$ is $m$%
\textbf{-divisible} if it can be written as $X=Y_{1}+Y_{2}+\cdots+Y_{m}$,
where $Y_{1},Y_{2},\ldots,Y_{m}$ are independent and identically distributed.
A\ random variable is \textbf{infinitely divisible} if it is $m$-divisible for
every $m$. \ We wish to study the allocation of $X=N$ cyclic points among a
fixed number $m$ of components, given a constrained random mapping. \ This
would be a matter of determining the inverse Laplace transform of the
$m^{\text{th}}$ root of%
\[
\mathcal{L}\left[  \sqrt{\frac{2}{\pi}}\exp\left(  \dfrac{-\xi^{2}}{2}\right)
\right]  =\exp\left(  \frac{\eta^{2}}{2}\right)  \operatorname*{erfc}\left(
\frac{\eta}{\sqrt{2}}\right)  .
\]
Pavlov's work \cite{Pav1-tcs0, Pav2-tcs0} is crucial here. \ We confront,
however, a surprising theoretical obstacle: the half-normal density is
provably \textit{not} infinitely divisible \cite{Ru-tcs0, Lu-tcs0, FH-tcs0}.
\ The independence requirement fails, in fact, beginning at $m=2$. \ Let us
offer a plausibility argument supporting this latter assertion.

On the one hand, if $Y_{2}=0$ is fixed, then the density of $Y_{1}=N$ clearly
approaches $\sqrt{2/\pi}$ as $\xi\rightarrow0^{+}$, i.e., it is bounded near
the origin.

On the other hand, if no condition is placed on $Y_{2}$, then the density of
$Y_{1}+Y_{2}=N$ is a convolution in the $\xi$-domain, which becomes
multiplication in the $\eta$-domain. \ Starting with \cite{Birn-tcs0,
Kom-tcs0, GU-tcs0}%
\[
\frac{\pi}{\sqrt{2\pi}+\pi\eta}<\sqrt{\frac{\pi}{2}}\exp\left(  \frac{\eta
^{2}}{2}\right)  \operatorname*{erfc}\left(  \frac{\eta}{\sqrt{2}}\right)
<\frac{\pi}{\sqrt{2\pi}+2\eta}%
\]
for all $\eta>0$, we deduce%
\[
\sqrt{\frac{1}{1+\sqrt{\frac{\pi}{2}}\,\eta}}<\exp\left(  \frac{\eta^{2}}%
{4}\right)  \sqrt{\operatorname*{erfc}\left(  \frac{\eta}{\sqrt{2}}\right)
}<\sqrt{\frac{1}{1+\sqrt{\frac{2}{\pi}}\,\eta}}%
\]
and upper/lower bounds are tight approximations of the center for small/large
values of $\eta$. \ No closed-form expression for $\mathcal{L}^{-1}[$center$]$
seems to be possible;\ $\mathcal{L}^{-1}[$lower bound$]$ and $\mathcal{L}%
^{-1}[$upper bound$]$ are%
\[%
\begin{array}
[c]{ccccc}%
\dfrac{2^{1/4}}{\pi^{3/4}\sqrt{\xi}}\exp\left(  -\sqrt{\dfrac{2}{\pi}}%
\,\xi\right)  &  & \text{and} &  & \dfrac{1}{(2\pi)^{1/4}\sqrt{\xi}}%
\exp\left(  -\sqrt{\dfrac{\pi}{2}}\,\xi\right)
\end{array}
\]
respectively. \ Both expressions approach infinity as $\xi\rightarrow0^{+}$,
tentatively implying that the density of $Y_{1}$ is unbounded near the origin.
This contrasts with the behavior described earlier, i.e., information about
$Y_{2}$ truly affects how $Y_{1}$ is distributed. \ Therefore $Y_{1}$ and
$Y_{2}$ must be dependent.

When we employed the word \textquotedblleft obstacle\textquotedblright%
\ before, it reflected our intention to study order statistics $Z_{1}%
=\min\{Y_{1},Y_{2}\}$ and $Z_{2}=\max\{Y_{1},Y_{2}\}$, with a goal of
understanding the allocation process (partioning $N$ cyclic points into two
cycles). \ If $Y_{1}$ and $Y_{2}$ were independent with common density
$\varphi(y)$, then the joint density of $Z_{1}\leq Z_{2}$ would simply be
$2\varphi(z_{1})\varphi(z_{2})$ for $z_{1}\leq z_{2}$. \ Dependency renders
the analysis more complicated.

As an aside, the Rayleigh density is also \textit{not} infinitely divisible
\cite{FH-tcs0}. The independence requirement again fails beginning at $m=2$.
\ We argue as before, but less formally. \ On the one hand, if $Y_{2}=0$ is
fixed, then the density of $Y_{1}$ approaches $0$ as $\xi\rightarrow0^{+}$.
\ On the other hand, if no condition is placed on $Y_{2}$, then we wish to
find the inverse Laplace transform of the square root of%
\[
\mathcal{L}\left[  \xi\exp\left(  \dfrac{-\xi^{2}}{2}\right)  \right]
=1-\sqrt{\frac{\pi}{2}}\,\eta\exp\left(  \frac{\eta^{2}}{2}\right)
\operatorname*{erfc}\left(  \frac{\eta}{\sqrt{2}}\right)  .
\]
A remarkably accurate approximation (with error less than $0.5\%$):%
\[
\sqrt{1-\sqrt{\frac{\pi}{2}}\,\eta\exp\left(  \frac{\eta^{2}}{2}\right)
\operatorname*{erfc}\left(  \frac{\eta}{\sqrt{2}}\right)  }\approx\exp\left(
\frac{\eta^{2}}{\pi}\right)  \operatorname*{erfc}\left(  \frac{\eta}{\sqrt
{\pi}}\right)
\]
defies easy explanation and yet provides a very helpful estimate:
\[
\mathcal{L}^{-1}\left[  \exp\left(  \frac{\eta^{2}}{\pi}\right)
\operatorname*{erfc}\left(  \frac{\eta}{\sqrt{\pi}}\right)  \right]
=\exp\left(  \frac{-\pi\,\xi^{2}}{4}\right)
\]
which approaches $1$ as $\xi\rightarrow0^{+}$. \ Since $0\neq1$, it follows
that $Y_{1}$ and $Y_{2}$ must be dependent.

\section{Addendum: Fallibility}

With the benefit of hindsight, we should have focused not on $\sigma(x)$, but
instead on
\[
\tilde{\sigma}(x)=\sqrt{x}\,\sigma(x)
\]
both here and in \cite{Fi4-tcs0}. \ The derivative of $\tilde{\sigma}(1/x)$ is
found as follows:%
\[
\frac{d}{dx}\tilde{\sigma}\left(  \frac{1}{x}\right)  =\frac{d}{dx}\left[
\frac{1}{x^{1/2}}\sigma\left(  \frac{1}{x}\right)  \right]  =-\frac
{1}{2x^{3/2}}\sigma\left(  \frac{1}{x}\right)  +\frac{1}{x^{1/2}}%
\sigma^{\prime}\left(  \frac{1}{x}\right)  \left(  -\frac{1}{x^{2}}\right)
\]
but%
\[
\sigma^{\prime}(y)=-\frac{1}{2y}\left(  \sigma(y)+\sigma(y-1)\right)
\]
therefore%
\begin{align*}
\frac{d}{dx}\tilde{\sigma}\left(  \frac{1}{x}\right)   &  =-\frac{1}{2x^{3/2}%
}\sigma\left(  \frac{1}{x}\right)  +\frac{1}{x^{1/2}}\left(  -\frac{x}%
{2}\right)  \left[  \sigma\left(  \frac{1}{x}\right)  +\sigma\left(  \frac
{1}{x}-1\right)  \right]  \left(  -\frac{1}{x^{2}}\right) \\
&  =-\frac{1}{2x^{3/2}}\sigma\left(  \frac{1}{x}\right)  +\frac{1}{2x^{3/2}%
}\left[  \sigma\left(  \frac{1}{x}\right)  +\sigma\left(  \frac{1}%
{x}-1\right)  \right]  =\frac{1}{2x^{3/2}}\sigma\left(  \frac{1-x}{x}\right)
,
\end{align*}
as was to be shown. \ Finally, $\rho(x)$ and $\sigma(x)$ are subsumed by the
general DDE \cite{Watt-tcs0, ABT1-tcs0, ABT2-tcs0}%
\[%
\begin{array}
[c]{ccc}%
x\,g^{\prime}(x)+(1-\theta)g(x)+\theta\,g(x-1)=0\text{ for }x>1, &  &
g(x)=x^{\theta-1}\text{ for }0<x\leq1
\end{array}
\]
where $\theta>0$ is fixed, and its associated Laplace transform is%
\[%
\begin{array}
[c]{ccc}%
\mathcal{L}\left[  g(\xi)\right]  =\dfrac{\exp\left(  -\theta\,E(\eta)\right)
}{\eta^{\theta}/\Gamma(\theta)}, &  & \eta\in\mathbb{C}\smallsetminus
(-\infty,0].
\end{array}
\]
It would be good someday to learn, from an interested reader, about possible
random mapping-theoretic applications of $g(x)$\ for select$\ \theta
\notin\{\frac{1}{2},1,\frac{3}{2},2,\frac{5}{2},\ldots\}$.

\section{Acknowledgements}

I am grateful to Ljuben Mutafchiev \cite{Muta-tcs0}, Lennart Bondesson
\cite{Bo1-tcs0, Bo2-tcs0} and Jeffrey Steif \cite{FSW-tcs0} for helpful
discussions. \ Vaclav Kotesovec provided relevant asymptotic expansions in
\cite{O1-tcs0}. \ The Mathematica routines NDSolve for DDEs \cite{Bel-tcs0}
and NInverseLaplaceTransform \cite{Val-tcs0} assisted in numerically
confirming many results. \ Interest in this subject has, for me, spanned many
years \cite{Fi5-tcs0, Fi6-tcs0}.


\begin{thebibliography}{99}                                                                                               %


\bibitem {Dick-tcs0}K. Dickman, On the frequency of numbers containing prime
factors of a certain relative magnitude, \textit{Ark. Mat. Astron. Fysik}, v.
22A (1930) n. 10, 1--14.

\bibitem {More-tcs0}P. Moree, Nicolaas Govert de Bruijn, the enchanter of
friable integers, \textit{Indag. Math.} 24 (2013) 774--801; arXiv:1212.1579; MR3124806.

\bibitem {Lag-tcs0}J. C. Lagarias, Euler's constant: Euler's work and modern
developments, \textit{Bull. Amer. Math. Soc.} 50 (2013) 527--628;
arXiv:1303.1856; MR3090422.

\bibitem {Watt-tcs0}G. A. Watterson, The stationary distribution of the
infinitely-many neutral alleles diffusion model, \textit{J. Appl. Probab.} 13
(1976) 639--651; 14 (1977) 897; MR0504014 and MR0504015.

\bibitem {ABT1-tcs0}R. Arratia, A. D. Barbour and S. Tavar\'{e}, Random
combinatorial structures and prime factorizations, \textit{Notices Amer. Math.
Soc.} 44 (1997) 903--910; MR1467654.

\bibitem {ABT2-tcs0}R. Arratia, A. D. Barbour and S. Tavar\'{e},
\textit{Logarithmic Combinatorial Structures: a Probabilistic, Approach},
Europ. Math. Society, 2003, pp. 21-24, 52, 87--89, 118; MR2032426.

\bibitem {Fi1-tcs0}S. R. Finch, Permute, Graph, Map, Derange, arXiv:2111.05720.

\bibitem {Fi2-tcs0}S. R. Finch, Rounds, Color, Parity, Squares, arXiv:2111.14487.

\bibitem {Fi3-tcs0}S. R. Finch, Second best, Third worst, Fourth in line, arXiv:2202.07621.

\bibitem {Fi4-tcs0}S. R. Finch, Joint probabilities within random
permutations, arXiv:2203.10826.

\bibitem {Step-tcs0}V. E. Stepanov, Limit distributions of certain
characteristics of random mappings (in Russian), \textit{Teor. Verojatnost. i
Primenen.} 14 (1969) 639--653; Engl. transl. in \textit{Theory Probab. Appl.
}14 (1969) 612--626; MR0278350.

\bibitem {OB-tcs0}F. Oberhettinger and L. Badii, \textit{Tables of Laplace
Transforms}, Springer-Verlag, 1973, pp. 172, 364, 367; MR0352889. \ 

\bibitem {Muta-tcs0}L. R. Mutafchiev, Large components and cycles in a random
mapping pattern, \textit{Random Graphs '87}, Proc. 1987 Pozna\'{n} conf., ed.
M. Karo\'{n}ski, J. Jaworski and A. Ruci\'{n}ski, Wiley, 1990, pp. 189--202; MR1094133.

\bibitem {Kolc-tcs0}V. F. Kolchin, \textit{Random Mappings}, Optimization
Software, 1986, pp. 35, 46--48, 85--88, 93--94, 153; MR0865130.

\bibitem {Pap-tcs0}A. Papoulis, \textit{Probability, Random Variables, and
Stochastic Processes}, McGraw-Hill, 1965, pp. 174--179; MR0176501.

\bibitem {RS-tcs0}H. Rubin and R. Sitgreaves, Probability distributions
related to random transformations of a finite set, Stanford Univ. tech.
report, 1954, http://statistics.stanford.edu/technical-reports/probability-distributions-related-random-transformations-finite-set.

\bibitem {Har-tcs0}B. Harris, Probability distributions related to random
mappings, \textit{Annals Math. Statist.} 31 (1960) 1045--1062; MR0119227.

\bibitem {PW-tcs0}P. W. Purdom and J. H. Williams, Cycle length in a random
function, \textit{Trans. Amer. Math. Soc.} 133 (1968) 547--551; MR0228032.

\bibitem {SL-tcs0}L. A. Shepp and S. P. Lloyd, Ordered cycle lengths in a
random permutation, \textit{Trans. Amer. Math. Soc.} 121 (1966) 340--357; MR0195117.

\bibitem {Grf-tcs0}R. C. Griffiths, On the distribution of allele frequencies
in a diffusion model, \textit{Theoret. Population Biol. }15 (1979) 140--158; MR0528914.

\bibitem {Shi-tcs0}T. Shi, Cycle lengths of $\theta$-biased random
permutations, B.S. thesis, Harvey Mudd College, 2014, http://scholarship.claremont.edu/hmc\_theses/65/.

\bibitem {Ren-tcs0}A. R\'{e}nyi, On connected graphs. I, \textit{Magyar Tud.
Akad. Mat. Kutat\'{o} Int. K\"{o}zl.} 4 (1959) 385--388; MR0126842.

\bibitem {O1-tcs0}N. J. A. Sloane, On-Line Encyclopedia of Integer Sequences,
A001865, A060281, A065456, and A273434.

\bibitem {Pav1-tcs0}Yu. L. Pavlov, Limit theorems for a characteristic of a
random mapping (in Russian), \textit{Teor. Veroyatnost. i Primenen.} 26 (1981)
841--847; Engl. transl. in \textit{Theory Probab. Appl. }26 (1982) 829--834; MR0636781.

\bibitem {Kru-tcs0}M. D. Kruskal, The expected number of components under a
random mapping function, \textit{Amer. Math. Monthly} 61 (1954) 392--397; MR0062973.

\bibitem {Rio-tcs0}J. Riordan, Enumeration of linear graphs for mappings of
finite sets, \textit{Annals Math. Statist.} 33 (1962) 178--185; MR0133250.

\bibitem {Pav2-tcs0}Yu. L. Pavlov, Random mappings with constraints on the
number of cycles (in Russian), \textit{Trudy Mat. Inst. Steklov.} 177 (1986)
122--132, 208; Engl. transl. in \textit{Proc. Steklov Inst. Math.} 1988, n. 4,
131--142; MR0840680.

\bibitem {De1-tcs0}J. M. DeLaurentis, Random functions and their small cycles,
\textit{Proc. 18}$^{\text{th}}$\textit{ Southeastern Conf. on Combinatorics,
Graph Theory and Computing}, Boca Raton, 1987, ed. F. Hoffman, R. C. Mullin,
R. G. Stanton and K. Brooks Reid, Congr. Numer. 58, Utilitas Math., 1987, pp.
295--302; MR0944710.

\bibitem {De2-tcs0}J. M. DeLaurentis, Components and cycles of a random
function, \textit{Advances in Cryptology -- CRYPTO '87}, Proc. 1987 Santa
Barbara conf., Lect. Notes in Comp. Sci. 293, Springer-Verlag, 1988, pp.
231--242; MR0956654.

\bibitem {RZ1-tcs0}D. Romero and F. Zertuche, The asymptotic number of
attractors in the random map model, \textit{J. Phys. A} 36 (2003) 3691--3700; MR1984723.

\bibitem {RZ2-tcs0}D. Romero and F. Zertuche, Grasping the connectivity of
random functional graphs, \textit{Studia Sci. Math. Hungar.} 42 (2005) 1--19; MR2128671.

\bibitem {MPQS-tcs0}R. S. V. Martins, D. Panario, C. Qureshi and E. Schmutz,
Periods of iterations of functions with restricted preimage sizes, \textit{ACM
Trans. Algorithms} 16 (2020) A. 30; arXiv:1701.09148; MR4120518.

\bibitem {Kup-tcs0}J. Kupka, The distribution and moments of the number of
components of a random function, \textit{J. Appl. Probab.} 27 (1990) 202--207; MR1039196.

\bibitem {Eng-tcs0}B. J. English, Factorial moments for random mappings by
means of indicator variables, \textit{J. Appl. Probab.} 30 (1993) 167--174; MR1206359.

\bibitem {O2-tcs0}N. J. A. Sloane, On-Line Encyclopedia of Integer Sequences, A201685.

\bibitem {FO-tcs0}P. Flajolet and A. M. Odlyzko, Random mapping statistics,
\textit{Advances in Cryptology - EUROCRYPT '89}, ed. J.-J. Quisquater and J.
Vandewalle, Lect. Notes in Comp. Sci. 434, Springer-Verlag, 1990, pp.
329--354; MR1083961.

\bibitem {Ru-tcs0}A. Ruegg, A characterization of certain infinitely divisible
laws, \textit{Annals Math. Statist.} 41 (1970) 1354--1356; MR0267626.

\bibitem {Lu-tcs0}E. Lukacs, \textit{Developments in Characteristic Function
Theory}, Macmillan Co., 1983, pp. 68--70; MR0810001.

\bibitem {FH-tcs0}F. W. Steutel and K. van Harn, \textit{Infinite Divisibility
of Probability Distributions on the Real Line}, Dekker, 2004, pp. 77, 113,
125--126, 411--412; MR.2011862.

\bibitem {Birn-tcs0}Z. W. Birnbaum, An inequality for Mill's ratio,
\textit{Annals Math. Statist.} 13 (1942) 245--246; MR0006640.

\bibitem {Kom-tcs0}Y. Komatu, Elementary inequalities for Mills' ratio,
\textit{Rep. Statist. Appl. Res. Un. Japan. Sci. Engrs.} 4 (1955) 69--70; MR0079844.

\bibitem {GU-tcs0}A. Gasull and F. Utzet, Approximating Mills ratio,
\textit{J. Math. Anal. Appl.} 420 (2014) 1832--1853; MR3240110.

\bibitem {Bo1-tcs0}L. Bondesson, On the infinite divisibility of the
half-Cauchy and other decreasing densities and probability functions on the
nonnegative line, \textit{Scand. Actuar. J.} (1987) 225--247; MR0943583.

\bibitem {Bo2-tcs0}L. Bondesson, \textit{Generalized Gamma Convolutions and
Related Classes of Distributions and Densities}, Lect. Notes in Stat. 76,
Springer-Verlag, 1992, p. 67; MR1224674.

\bibitem {FSW-tcs0}B. G. Franz\'{e}n, J. E. Steif and J. W\"{a}stlund, Where
to stand when playing darts? \textit{ALEA Latin Amer. J. Probab. Math.
Statist.} 18 (2021) 1561--1583; arXiv:2010.00566; MR4291033.

\bibitem {Bel-tcs0}Y. M. Beltukov, Tuning NDSolve options to accurately solve
DDEs, http://mathematica.stackexchange.com/questions/59431/simpler-code-evaluating-dickmans-function.

\bibitem {Val-tcs0}P. Valko, NInverseLaplaceTransform numerical approximation
package, http://resources.wolframcloud.com/FunctionRepository/resources/NInverseLaplaceTransform/.

\bibitem {Fi5-tcs0}S. R. Finch, Golomb-Dickman constant, \textit{Mathematical
Constants}, Cambridge Univ. Press, 2003, pp. 284--292; MR2003519.

\bibitem {Fi6-tcs0}S. R. Finch, Extreme prime factors, \textit{Mathematical
Constants II}, Cambridge Univ. Press, 2019, pp. 171--172; MR3887550.%

\begin{tabular}
[c]{lll}
& Steven Finch & \\
& MIT Sloan School of Management & \\
& Cambridge, MA, USA & \\
& \textit{steven\_finch@harvard.edu} &
\end{tabular}

\end{thebibliography}
\end{document}